\documentclass[11pt,leqno]{amsart}
\usepackage{amsmath,amsfonts,amssymb,amscd,amsthm,amsbsy,epsf}
\newtheorem{thm}{Theorem}
\newtheorem{lem}[thm]{Lemma}

\theoremstyle{remark}

\newcommand{\R}{\mathbb R}
\newcommand{\eps}{\varepsilon}

\newcommand{\lap} {\triangle}
\DeclareMathOperator{\tr}{tr}
\begin{document}
\title{On the Evans-Krylov theorem}
\author{Luis  Caffarelli \and Luis Silvestre}
\address{Luis Caffarelli\\
Department of Mathematics\\
University of Texas at Austin\\
1 University Station -- C1200\\
Austin, TX~78712-0257}
\email{caffarel@math.utexas.edu}
\address{Luis Silvestre\\
Mathematics Department\\
University of Chicago\\
Chicago, IL 60637}
\email{luis@math.uchicago.edu}

\maketitle
\pagestyle{plain}

The Evans-Krylov theorem consists of the a~priori estimate:

\begin{thm} \label{t:main}
Smooth solutions, $u$, of a uniformly elliptic, fully non-linear convex 
equation $F(D^2 u) =0$ in the unit ball $B_1$, of $\R^n$ have a 
$C^{2,\alpha}$ interior a~priori estimate 
$$\|u\|_{C^{2,\alpha} (B_{1/2})} 
\le C\| u\|_{C^{1,1} (B_1)}$$
with the constant $C$ depending only on the ellipticity 
of $F$.
\end{thm}

The importance of the Evans-Krylov theorem is that it allows us to solve 
the Dirichlet problem for fully nonlinear equations by the method of 
continuity (rendering classical solutions). 

This theorem was proved independently by N. Krylov \cite{MR661144} and L. C. Evans \cite{MR649348}.
In this note, motivated by our work on integral fully nonlinear equations \cite{CS3}, we 
provide a more direct presentation of their proof (although the underlying 
key ideas are the same).

We recall the two opposite components in Krylov-Safonov Harnack inequality, the proof of which can be found in \cite{caffarelli1995fne} (Theorem 4.8).

\begin{enumerate}
\item[a)] (The weak $L^\eps$ estimate)
If $v$ is a non-negative supersolution of 
$$a_{ij} (x) D_{ij} v \le 0$$
in $B_1$, with $\lambda I \le a_{ij} \le \Lambda I$ then 
$$|\{ v>t \inf_{B_{1/2}} v \} \cap B_{1/4} | \le C(\lambda,\Lambda) t^{-\eps}$$

\item[b)] (the oscillation lemma)
If $v$ is a subsolution of $a_{ij} (x) D_{ij} v \ge 0$ in $B_1$ and $v\le 1\ ,$
then
$$\sup_{B_{1/2}} v \le C(\lambda,\Lambda) |\{ v > 0\} \cap B_{3/4} |$$ 
\end{enumerate}

In case of harmonic functions, these are just consequences of the mean value theorem.

We also recall that convexity of $F$ as a function of $D^2u$, implies that 
any pure second derivative, $u_{\sigma\sigma}$, of $u$ and thus any linear 
combination
$$\ell (x) = \sum_j u_{\sigma_j\sigma_j} (x)$$
is a supersolution of the linearized operator 
$$a_{ij} (x) D_{ij} \ell (x) \le 0$$
($a_{ij} (x) = F_{ij} (D^2 u(x)$).

Finally, the uniform ellipticity of $F$ implies that for any two points
$x_1, x_2$ in $B_1$,
\begin{equation} \label{e:comparable}
\tr [ D^2 u (x_2) - D^2 u(x_1)]^+  \approx \tr [D^2 u(x_2 ) - D^2 u (x_1)]^-
\end{equation}

At this point, we define for any subspace $V$
$$w(x,V) = \Delta_V u (x) - \Delta_V u (0)$$
($\Delta_V u(x)$ is the Laplacian of $u$ at the point $x$ when restricted to the affine variety $x+V$).

Note that for each fixed $V$, $w$ is an $\ell(x)$ as above and satisfies the $L^\eps$ estimate. Also, note that the \emph{positive} and \emph{negative} part of the laplacian can be expressed as
\begin{align*}
\max_V w(x,V) &=  \tr [ D^2 u (x) - D^2 u(0)]^+, \\
\min_V w(x,V) &=  -\tr [ D^2 u (x) - D^2 u(0)]^-.
\end{align*}
By rescaling dyadically and iterating it is enough to prove the 
following lemma:

\begin{lem}
There exists a $\theta > 0$, $\theta = \theta (\lambda, \Lambda)$, such that if 
for all $V$, for all $x$ in $B_1$, 
$$w(x,V) \ge -1$$
Then for all $V$, for all $x$ in $B_{1/2}$,
$$w(x,V) \ge -1 + \theta .$$
\end{lem}

Indeed, this will imply by iteration, that the laplacian is H\"older continuous. Noew we prove the lemma.

\begin{proof}
Assume that $w(x_0,V_0) \le -1 +\theta$ for some $V_0$ and $x_0$ in $B_{1/2}$ 
($\theta$, small, to be chosen). We will then find a contradiction. Since $w(\cdot,V) +1$ 
is a nonnegative supersolution the $L^\eps$ lemma applies and 
$$w(x,V) + 1 \le \theta^{1/2}$$
in a set $\Omega$ that covers \emph{almost all} of $B_{1/4}$, i.e., 
$$|B_{1/4} \setminus \Omega | \le C\theta^{\eps/2}$$

We notice that in $\Omega$, $1-\theta^{1/2} \leq -w(x,V) \leq \tr [ D^2 u (x) - D^2 u(0)]^- \leq 1$. On the other hand, we know that \[w(x,V) + w(x,V^\perp) = \lap u(x) - \lap u(0) = \tr [ D^2 u (x) - D^2 u(0)]^+ - \tr [ D^2 u (x) - D^2 u(0)]^-.\]
Thus, we also have $0 \leq \tr [ D^2 u (x) - D^2 u(0)]^+ - w(x,V^\perp) \leq \theta^{1/2}$ for $x \in \Omega$. Moreover, for $\theta$ small, by \eqref{e:comparable},
\[ -w(x,V) \approx \tr [ D^2 u (x) - D^2 u(0)]^- \approx \tr [ D^2 u (x) - D^2 u(0)]^+ \approx w(x,V^\perp).\]
Thus, there is a constant $c(\lambda,\Lambda)>0$ such that $w(x,V^\perp) \geq c(\lambda,\Lambda)$ in $\Omega$. We now examine the function $v = (c(\lambda,\Lambda) - w(x,V^\perp))^+$ in $B_{1/4}$, for which the oscillation lemma applies and satisfies
\begin{itemize}
\item[a)] $0\le v\le 2$
\item[b)] $v(0) = c(\lambda, \Lambda)$
\item[c)] $v = 0$ in $\Omega$.
\end{itemize}
For $\theta$ small (i.e., for $\Omega$ \emph{almost all} of $B_{1/4}$) this 
contradicts the oscillation lemma since $c(\lambda,\Lambda)$ is a fixed positive constant for $\theta$ small. This completes the proof.
\end{proof}

\bibliographystyle{plain}
\bibliography{evansKrylov.bib}
\end{document}